\documentclass[review, 3p,10pt, twocolumn]{elsarticle}

\usepackage{amssymb}
 \usepackage{amsmath}
 \usepackage{hyperref}

 \usepackage[ruled,vlined,english,linesnumbered]{algorithm2e}

\newtheorem{prop}{Proposition} % [section]
\newtheorem{lemma}[prop]{Lemma}
\newtheorem{theorem}[prop]{Theorem}

\usepackage{color}

\definecolor{mygray}{gray}{0.6}

\begin{document}

\begin{frontmatter}

%% Title, authors and addresses

%% use the tnoteref command within \title for footnotes;
%% use the tnotetext command for theassociated footnote;
%% use the fnref command within \author or \address for footnotes;
%% use the fntext command for theassociated footnote;
%% use the corref command within \author for corresponding author footnotes;
%% use the cortext command for theassociated footnote;
%% use the ead command for the email address,
%% and the form \ead[url] for the home page:
%% \title{Title\tnoteref{label1}}
%% \tnotetext[label1]{}
%% \author{Name\corref{cor1}\fnref{label2}}
%% \ead{email address}
%% \ead[url]{home page}
%% \fntext[label2]{}
%% \cortext[cor1]{}
%% \address{Address\fnref{label3}}
%% \fntext[label3]{}

%\title{\rev{Cactus metrics and their optimal realizations}}
\title{Recognizing and realizing cactus metrics}

%% use optional labels to link authors explicitly to addresses:
%% \author[label1,label2]{}
%% \address[label1]{}
%% \address[label2]{}

%\author{}

%\address{}

\author[tok,jst]{Momoko Hayamizu}
\ead{hayamizu@ism.ac.jp}

\author[uea]{Katharina T. Huber}
\ead{k.huber@uea.ac.uk}

\author[uea]{Vincent Moulton\corref{cor1}}
\ead{v.moulton@uea.ac.uk}

\author[tu]{Yukihiro Murakami}
\ead{y.murakami@tudelft.nl}

\address[tok]{The Institute of Statistical Mathematics, 190-‐8562, Tachikawa, Tokyo, Japan}

\address[jst]{JST PRESTO, 190-‐8562, Tachikawa, Tokyo, Japan}

\address[uea]{School of Computing Sciences, University of East Anglia, NR4 7TJ, Norwich, United Kingdom}

\address[tu]{Delft Institute of Applied Mathematics, Delft University of Technology, Van Mourik Broekmanweg 6, 2628 XE, Delft,The Netherlands}

\cortext[cor1]{Corresponding author}

\begin{abstract}
The problem of realizing finite metric spaces in terms of weighted graphs 
has many applications. For example, the mathematical and computational properties 
of metrics that can be realized by trees have been well-studied 
and such research has laid the foundation of the reconstruction of phylogenetic trees from evolutionary distances. 
However, as trees may be too restrictive to accurately represent real-world data or phenomena, it 
is important to understand the relationship between more general graphs and distances. In this paper, 
we introduce a new type of metric called a cactus metric, that is, a metric that can 
be realized by a cactus graph. We show that, just as with tree metrics, a cactus 
metric has a unique optimal realization. In addition, we describe an
%polynomial time 
algorithm that can recognize whether or not a metric is a cactus 
metric and, if so, compute its optimal realization in $O(n^3)$ time, where $n$
is the number of points in the space. 
\end{abstract}

\begin{keyword}
Cactus metric \sep $X$-cactus \sep Metric realization \sep Optimal realization \sep Phylogenetic network  
%% keywords here, in the form: keyword \sep keyword

%% PACS codes here, in the form: \PACS code \sep code

%% MSC codes here, in the form: \MSC code \sep code
%% or \MSC[2008] code \sep code (2000 is the default)

\MSC 05C05 \sep  05C12 \sep 92B10  
\end{keyword}

\end{frontmatter}

%%%%%%%%%%%%%%%%%%%%%%%%%%%%%%%%%%%%%%%%%%%%%%%%%%%%%%%%%%%%%%%%%%%%%%%%%%
\section{Introduction}
\label{section:introduction}
%%%%%%%%%%%%%%%%%%%%%%%%%%%%%%%%%%%%%%%%%%%%%%%%%%%%%%%%%%%%%%%%%%%%%%%%%%

The metric realization problem, which is the problem of representing a finite 
metric space by a weighted graph, has many applications, most notably in
the reconstruction of evolutionary trees. Although any finite metric space 
can be realized by a weighted complete graph, there can be different graphs 
that induce the same metric. In \cite{hak}, Hakimi and Yau first 
considered “optimal” realizations of finite metric spaces, 
%that is, those that are shortest among all possible realizations. 
which are realizations of least total weight.  Although every finite metric space 
has an optimal realization \cite{dress, imrich84}, the problem of finding an 
optimal realization is NP-hard in general \cite{alt, wink} and the optimal 
solution is not necessarily unique \cite{alt, dress}. 

A well-known special case of optimal realizations is provided by tree metrics, namely, those 
metrics that can be realized by some edge-weighted tree. For any 
tree metric on a finite set $X$, its optimal realization is an $X$-tree (\textit{i.e.}, a tree 
in which some vertices are  labeled by $X$) and is uniquely determined \cite{hak}. 
In addition, there exist optimal polynomial-time algorithms for computing 
the tree realization from a tree metric \cite{band,bat,culb}. However, not much is 
known about the properties of optimal realizations of metrics induced 
by graphs that are more general than trees. Developing our understanding 
in this direction could be useful, as trees can sometimes be too restrictive for 
realizing metrics arising in real-world applications \cite{huson}.

\begin{figure}[htbp]
\centering
\includegraphics[width=0.35\textwidth]{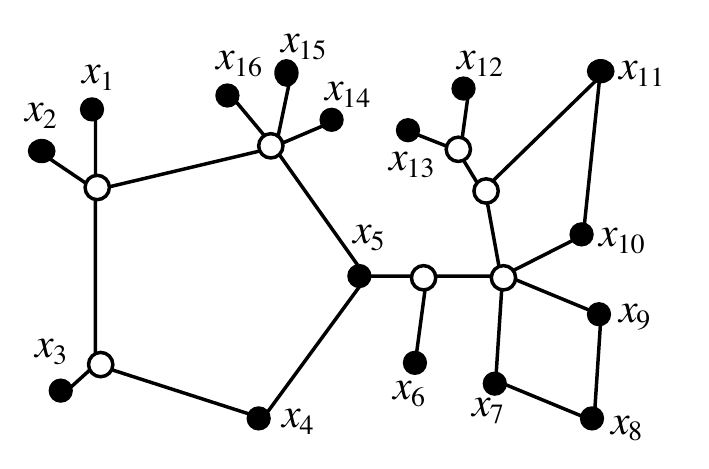}
\caption{An example of an $X$-cactus with a label-set $X=\{x_1,\dots,x_{16}\}$, where 
	the weight of each edge is proportional to its length. The vertices labeled by an 
	element of $X$ are shown in black. The white circles are vertices that are not in $X$.  
\label{fig:X-cactus}}
\end{figure}

In this paper, we generalize the concept of a tree metric  by  introducing a 
new type of metric called a ``cactus metric\footnote{This concept was 
	first introduced in \cite{hay}.}'' which can be realized by an 
edge-weighted ``$X$-cactus'', where a cactus is a connected graph 
in which each edge belongs to at most one cycle.
An example of an $X$-cactus is presented in Figure \ref{fig:X-cactus}. Note that 
cacti have some nice properties in common with trees. For instance,  every 
cactus is planar and the number of vertices in an $X$-cactus 
is $O(|X|)$ as with $X$-trees, which 
means that cactus metrics are easy to visualize. In particular, they provide a special case of an open 
problem in discrete geometry from  Matou{\v{s}}ek \cite{mat}. Besides these 
observations, in this paper we prove that, just as with tree metrics, any cactus 
metric has a unique optimal realization. We also describe a polynomial time 
algorithm for deciding whether or not an arbitrary metric is a cactus metric, which 
also computes its optimal realization in case it is.

\section{Preliminaries}\label{sec:prelim}

A {\em metric} on a set $S$ is defined to be a function $d: S \times S \to {\mathbb R}_{\ge 0}$ with 
the property that $d$ equals zero if and only if the two elements in $S$ are identical, is symmetric,  
and satisfies the triangle inequality.  

All graphs considered here are  finite, connected, simple, undirected graphs in which the edges 
have positive weights. % and some vertices are labeled using the label-set $X$.
For any graph $G$, $V(G)$ and $E(G)$ represent the vertex-set and edge-set of $G$, respectively. 
For any vertex $v$ of a graph $G$, the number of edges of $G$ that have $v$ as an endvertex is denoted by $deg(v)$. 
For any graph $G$ and any subset $S$ of $V(G)$, we let $d_G$ denote the metric 
on $S$ induced by taking shortest paths in $G$ between elements in $S$.

Throughout this paper, we use the symbol $X$ to 
represent a  finite set with $|X|\ge 2$, which is sometimes called a {\em label-set}. 
For any metric $d$ on $X$, a {\em realization} of $(X,d)$ is a graph $G$ 
such that $X$ is a subset of $V(G)$ and $d(x,y)=d_G(x,y)$ holds for each $x,y \in X$, where 
we shall always assume that each vertex $v$ of $G$ with $deg(v)\le 2$ has a label in $X$ \cite{imrich84}. 
A realization is {\em minimal} if the  removal of an arbitrary edge of $G$ yields a graph 
that does not realize $d$. It is {\em optimal} if the sum of 
its edge weights is minimum over all possible realizations (note that
optimal realizations are minimal but the converse does not hold). Any finite 
metric space has at least one optimal realization \cite[Theorem 2.2]{imrich84}.

We now state a theorem concerning optimal realizations which will be useful in our proofs.
For a graph $G$, each maximal biconnected subgraph of $G$ is called a {\em block} of $G$ 
and each vertex of $G$ shared by two or more blocks of $G$ is called a {\em cutvertex} of $G$. 
Notice that if a graph consists of a single block, then it has no cutvertex. 

\begin{theorem}[\cite{imrich84}, Theorem 5.9]~\label{optimalblock} 
	Let $G$ be a minimal realization of a finite metric space $(X, d)$, let 
	$G_1,\dots,G_k$ be the blocks of $G$, let $M_i$ be the union of the vertices of $X$ in $G_i$ together with the 
	cutvertices of $G$ in $G_i$, and let $d_i$ be the metric induced by $G$ on $M_i$. 
	Then,  if  every $G_i$ is an optimal realization of $(M_i, d_i)$, then $G$ is also optimal. If every $G_i$, besides being optimal, is also unique, then $G$ is optimal and unique too. 
\end{theorem}

We now turn to two special classes of metrics, that is, tree metrics and cyclelike metrics. 
A metric $d$ on $X$ is called a {\em tree metric}  if there exists an $X$-tree 
that realizes $(X, d)$, where an {\em $X$-tree} is a tree  $T$ with the property 
that each vertex $v$ of $T$ with $deg(v) \le 2$ is contained in $X$ \cite{ss13}. 

\begin{theorem}[\cite{hak}]\label{uniquetree}
If $d$ is a tree metric on a finite set $X$, then there exists an 
$X$-tree that is a unique optimal realization of $(X, d)$. 
\end{theorem}

Given a metric $d$ on $X$ with $|X|\ge 4$, we say that $d$ is {\em cyclelike} 
if there is a minimal realization for $d$ that is a cycle. 
This type of metric was discussed in  \textit{e.g.}, \cite{rubei18,imrich84, unicyclic}.
The following result will also be useful.

\begin{theorem}[\cite{imrich84}, Theorem 4.4]\label{optimalcycle}
	Suppose $d$ is a cyclelike metric on a finite set $X$ and a cycle $C$  is a minimal realization of $(X, d)$  with $V(C)=X = \{v_1,v_2,\dots,v_m\}$, $m\ge 4$,  
	and $E(C)=\{\{v_i,v_{i+1}\} \,:\, 1 \le i \le m\}$, where the indices are taken modulo $m$. Then, $C$ is an optimal realization of $(X, d)$ if and only if 
		$$
		d(v_{i-1},v_i) + d(v_i,v_{i+1}) = d(v_{i-1},v_{i+1})
		$$
		holds for all $i$. In this case, $C$ is the unique optimal realization of $(X, d)$.
\end{theorem}

\section{The uniqueness of optimal realizations of cactus metrics}

As mentioned above a {\em cactus} is a connected graph in which each edge belongs to at most one cycle.
We define an {\em $X$-cactus} to be a cactus $G$ with the property 
that each vertex $v$ of $G$ with $deg(v) \le 2$ is contained in $X$ (see Figure \ref{fig:X-cactus}).
Note that the maximum number of cycles in an $X$-cactus 
is $|X|-2$ (which can be proved by induction on $|X|$).
In addition, we say that a metric $d$ on a finite set $X$ is a {\em cactus metric}  if there exists 
an {\em edge-weighted} $X$-cactus that realizes $(X, d)$. 

Given an edge-weighted cycle $C= v_1,\dots,v_m$ that is a realization of its 
	corresponding metric $d_C$, we call a vertex $v_i \in V(C)$ {\em slack}
	if $d(v_{i-1},v_i) + d(v_i,v_{i+1}) > d(v_{i-1},v_{i+1})$. 
The following lemma is a direct consequence of Theorem \ref{optimalcycle}. 

\begin{lemma}\label{noslack}
	Under the premise of Theorem \ref{optimalcycle}, $C$ is an optimal realization of $(X, d)$ if and only if $C$ has no slack vertex.
\end{lemma}

We now use the lemma to prove the following generalization of Theorem~\ref{uniquetree}, using 
	the concept of ``compactification'' \cite{hak,unicyclic,compactification}. 
 
\begin{theorem}\label{unique}
If $d$ is a cactus metric on a finite set $X$, then there exists an $X$-cactus that is a unique optimal realization of $(X, d)$. 
\end{theorem}
\noindent{\em Proof: }
Let $G$ be an $X$-cactus that is a minimal realization of $(X, d)$.
Without loss of generality, we assume that each 
cycle of $G$ has at least four vertices (since we can
always replace a 3-cycle with a tree in such a way that the 
obtained graph is a realization). 
If there is no cycle in $G$ containing a slack vertex, then the 
assertion immediately follows from Theorems~\ref{optimalblock}, \ref{optimalcycle}  and Lemma~\ref{noslack}. 

\begin{figure*}[htbp]
\centering
\includegraphics[width=0.7\textwidth]{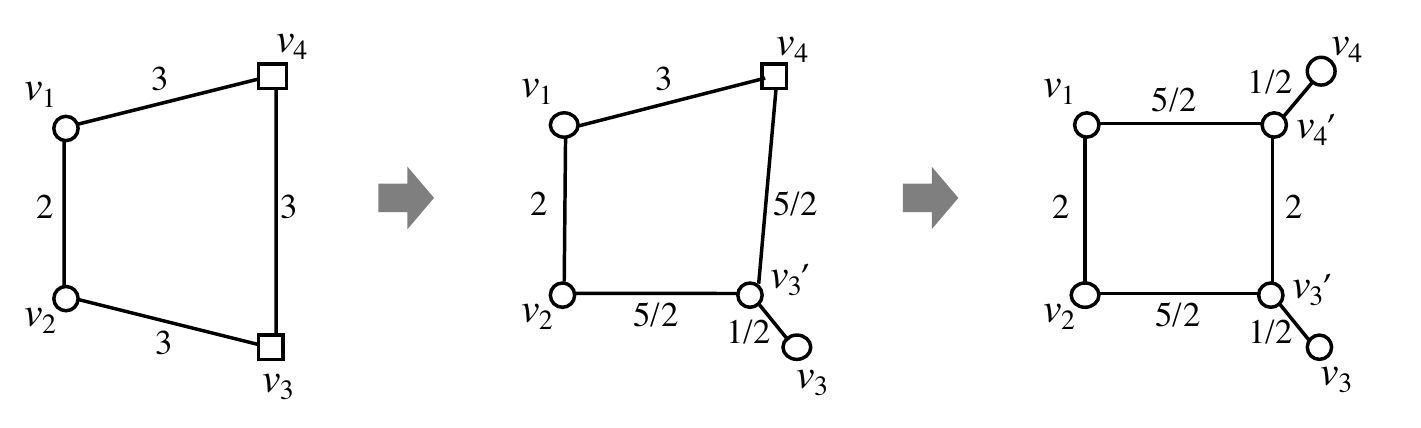}
\caption{An illustration of compactification that is described in the proof of Theorem \ref{unique}, where 
	we highlight each slack vertex by a square. Compactification of $v_3$ in the left graph yields 
	the graph in the middle panel, which still contains a slack vertex $v_4$. If we further 
	apply the same operation to $v_4$, then we obtain the graph on the right which has no slack vertex.
\label{fig:compactification}}
\end{figure*}

So, assume that there is a cycle $C=v_1,\dots,v_m$ in $G$ that  
has  consecutive edges $\{v_{i-1}, v_i\}, \{v_i, v_{i+1}\}$ with 
$\Delta_i := \{d_G(v_{i-1},v_i) + d_G(v_i,v_{i+1}) - d_G(v_{i-1},v_{i+1})\}/2>0$.
As we will now explain, we apply a ``compactification'' operation
to the slack vertex $v_i$ (see also Figure \ref{fig:compactification}).    
For notational convenience, let $\Delta_{i-1}:=\{d_G(v_{i-1},v_i) + d_G(v_{i-1},v_{i+1}) - d_G(v_i,v_{i+1})\}/2$ 
and $\Delta_{i+1}:=\{d_G(v_{i+1},v_i) + d_G(v_{i-1},v_{i+1}) - d_G(v_i,v_{i-1})\}/2$. 
Compactification of $v_i$ refers to converting $G$ into the graph $G^\prime$  
with $V(G^\prime):=V(G)\cup \{v_i^\prime\}$ and 
$E(G^\prime):=(E(G)\setminus \{\{v_{i-1}, v_i\}, \{v_i, v_{i+1}\}\})\cup \{\{v_{i-1}, v_i^\prime\}, \{v_i, v_i^\prime\}, \{v_{i+1}, v_i^\prime\}\}$, 
where  for each $j\in \{i-1, i, i+1\}$,  the edge $\{v_j, v_i^\prime\}$ has weight $\Delta_j$. 
As can be easily verified,  $G^\prime$ is an $X$-cactus that is 
a minimal realization of $(X,d)$ with a strictly smaller number of slack vertices 
than $G$. Thus, as $|V(G)|$ is finite, by applying the same operation  
repeatedly and suppressing all unlabeled vertices of degree two (if any arise), we 
will eventually obtain an $X$-cactus that realizes $(X,d)$ 
without a slack vertex, which must be the unique optimal realization of $(X, d)$. 
\qed

It is interesting to see that for cactus metrics, we do not need to perform too many ``compactifications''
for each cycle in the above proof in light of the following observation.

\begin{prop}\label{atmosttwo}
If the premise of Theorem~\ref{optimalcycle} holds, then 
$C$ has at most two  slack vertices. In the case when there exist precisely two slack vertices, they are adjacent in $C$.  
\end{prop}
\noindent{Proof:}
Let $V(C)=\{v_1,\dots,v_m\}$ as in Theorem~\ref{optimalcycle}. Suppose $C$ 
has at least two slack vertices and assume that $v_i$ is a slack vertex, in other words, that
$d(v_{i-1},v_i) + d(v_i,v_{i+1}) > d(v_{i-1},v_{i+1})$ holds. 
As the path in $C$ from $v_{i-1}$ to $v_{i+1}$ that does not contain $v_i$ is the shortest path 
between $v_{i-1}$ and $v_{i+1}$, it follows that any $v\in V(C)\setminus \{v_{i-1},v_i,v_{i+1}\}$ is not slack. 
Now, suppose $v_{i-1}$ is a slack vertex. Then 
using a similar argument 
%as for $v_i$ 
by considering the shortest path between $v_{i-2}$ and $v_i$,
it follows that $v_{i+1}$ is not slack. So the only slack vertices 
are $v_i$ and $v_{i-1}$. The same argument applies to the case when $v_{i+1}$ is a slack vertex. 
\qed\\

\section{A polynomial time algorithm for finding the optimal cactus realization}

In this section we describe an algorithm, which for a metric $d$ on $X$,  
produces the unique optimal
 realization for $d$ that is an $X$-cactus or a message that there is no such realization in $O(|X|^3)$ time. 
This should be compared to tree metrics for which the same process
can be carried out in $O(|X|^2)$ time \cite{bat,culb}.

We begin by considering cyclelike metrics. Note that 
the characterization given in Theorem~\ref{optimalcycle} for when a realization of a cyclelike metric is optimal is not
sufficient to characterize cyclelike metrics, as pointed out in \cite{unicyclic}. Even so we have the following result
(which is related to 
Theorem~4.1 in \cite{rubei18}):

\begin{lemma}\label{cyclecheck}
    Given a metric $d$ on $X$, we can determine if there is an edge-weighted cycle $C$ that is an 
    	optimal realization of $(X, d)$ and, 	if so, compute $C$ in $O(|X|^2)$ time.  
\end{lemma}
\noindent{\em Proof: }
We describe an algorithm that takes an arbitrary metric $d$ on
$X$ as input, which in case $d$ has an optimal realization that is a cycle  
computes this cycle, and stops if this is not the case:
		
1) Start by finding a pair $\{v_0, v_1\}$  of distinct elements in $X$ such that $d(v_0, v_1)\le d(p, q)$ 
holds for any $\{p, q\}\in  \binom{X}{2}\setminus \{\{v_0, v_1\}\}$, and then set 
$e_1:=\{v_0, v_1\}$ and $w_1:=d(v_0, v_1)$. 
2) For each $j\in \{2, \dots, |X|-1 \}$, find all vertices $x\in X\setminus \{v_0, \dots, v_{j-1}\}$ 
with $d(v_{j-2}, v_{j-1})+d(v_{j-1}, x)=d(v_{j-2}, x)$.
Among these vertices, we let $v_j$ be the unique vertex $x$ that minimizes $d(v_{j-1},x)$.
If such a vertex does not exist, or if such a vertex does exist 
but it is not unique, then stop; else set $e_j := \{v_{j-1}, v_j\}$ and $w_j:=d(v_{j-1}, v_j)$. 
3) Set $e_{|X|} := \{v_{|X|-1}, v_0\}$ and $w_{|X|}:=d(v_{|X|-1}, v_0)$. 
4) Check if the cycle $C$ defined by $V(C):=X$ and $E(C):=\{e_1,\dots, e_{|X|}\}$ together 
with the  weight $w_j$ of each edge $e_j\in E(C)$ is a minimal realization of $(X, d)$. If not then stop, else
output the weighted cycle $C$. 

If this algorithm returns a cycle  $C$ that realizes $(X, d)$, then $C$ satisfies 
the equation in Theorem \ref{optimalcycle} and so $C$ is the optimal 
realization of $(X, d)$. 
Conversely, if there is a cycle $C$ that is an optimal realization of $(X, d)$, 
then  $C$ is unique. In this case, the above algorithm correctly constructs $C$ as follows. 
The algorithm initializes by finding two vertices of $X$ that are closest together.
Since an optimal realization that is a cycle is minimal, it must be the case 
that these two vertices are connected by an edge.
In Step 2, the algorithm iteratively extends the existing path by seeking 
for the neighbour of $v_{j-1}$, which is one of the endvertices of the path.
Observe that the two conditions in Step 2 uniquely determine this neighbour: the 
first condition ensures that a shortest path between $v_{j-2}$ and $v_j$ contains $v_{j-1}$; the 
second condition correctly identifies the neighbour of $v_{j-1}$ by making 
sure that the distance between it and $v_{j-1}$ is shortest.
In Step 3, we join the two endvertices of the path by an edge to form the cycle $C$. 
Note that in this step, we run the risk of making a realization of $(X,d)$ that is a 
path into a realization of $(X,d)$ that is a cycle that is not minimal.
Due to this, and also to ensure we have the correct solution, we check 
that the cycle is a minimal realization of $(X,d)$ in Step 4.

To give the running time of the algorithm, observe that Step 1 
takes $O(|X|^2)$ time as we search for a minimum element from a set of size $\binom{|X|}{2}$.
In Step 2, we iterate over a `for loop' at most $|X|$ times. Within the `for loop' 
we iterate over at most $|X|$ elements to find the vertices that satisfy the first condition.
Then, we iterate over those vertices to find a minimum element from at most $|X|$ elements.
Hence, each `for loop' takes $O(|X|)$ time; it follows then that Step 2 takes $O(|X|^2)$ time.
Step 3 takes constant time, as we simply add a weighted edge to the graph.
Since one can obtain the metric induced by a cycle in at most $O(|X|^2)$ time, 
Step 4 can be performed in at most $O(|X|^2)$ time.
As each step of the algorithm can be done in $O(|X|^2)$ time, the whole 
algorithm requires $O(|X|^2)$ time. 
\qed\\

\begin{theorem} \label{poly}
	Given a metric $d$ on $X$, we can determine if $d$ is a cactus metric 
	and if so construct its optimal  realization in $O(|X|^3)$ time.  
\end{theorem}
\noindent{\em Proof: }
In  \cite[Algorithm 2]{HVcutpoint} Hertz and Varone 
give a polynomial time algorithm
for decomposing an arbitrary metric space $(X, d)$ into finite 
metric spaces $(M_i, d_i)$, $1 \le i \le k$, with $|M_i| \le |X|$, 
such that any optimal realization of $(M_i, d_i)$ 
must consist of a single block, and such that an optimal 
realization for $d$ can be constructed by piecing together the 
optimal realizations for the $(M_i, d_i)$. They
also observe \cite[p.174]{HVcutpoint} that this decomposition
can be computed in  $O(|X|^3)$ time using results in \cite{dress10} (see
also \cite[p.160]{dress10}). 
In addition, by the arguments in \cite[Lemma 3.1]{dress10}, it follows 
that $k$ is $O(|X|)$.

Assume that we have decomposed $(X, d)$ into $\{(M_i, d_i)\}_{i\in \{1,\dots,k\}}$ by 
using the aforementioned preprocessing algorithm. In case $|M_i|=2$, its 
optimal realization is obviously a tree. Recalling the argument 
in the proof of Theorem~\ref{unique}, we  know that $|M_i|\neq 3$  holds  
for each $i\in \{1,\dots,k\}$.
For each $(M_i, d_i)$ with $|M_i| \ge 4$, by using the algorithm in 
Lemma~\ref{cyclecheck}, we can check if $(M_i, d_i)$ has 
an optimal realization that is a cycle or not, and if so 
construct the cycle in $O(|M_i|^2)$ time (and hence  $O(|X|^2)$ time suffices). 
If there is some $i\in \{1,\dots,k\}$ such that $|M_i| \ge 4$ and $(M_i, d_i)$ does not 
have an optimal realization that is a cycle, then $d$ is not a cactus 
metric, else $d$ is a cactus metric, and we can construct 
the cactus by piecing together the optimal realizations 
for the $(M_i, d_i)$.  Using the aforementioned fact 
that $k$ is $O(|X|)$, we 
conclude that the overall time complexity is $O(|X|^3)$.
\qed\\

\section{Discussion and future work}

It may be worth investigating as to whether there is a more direct and efficient 
algorithm than the one given in Theorem~\ref{poly} for 
recognizing and/or realizing  cactus metrics that use
structural properties of cactus graphs.
More generally, we
could investigate optimal realizations for metrics that
can be realized by graphs $G$ in which every block $G_i=(V_i, E_i)$
satisfies $|E_i|-|V_i|+1 \le k$, and such that every vertex in $G$ 
with degree at most 2 is contained in $X$.
Here, we note that 
in case $k=0$, $G$ is an $X$-tree, and in case $k=1$, $G$ is an $X$-cactus.
However, even in case $k=2$, there may be infinitely many
optimal realizations (e.g. the metric given in \cite[Fig. 15]{alt}). So it might
be interesting to first understand for $k\ge 2$  
which of these metrics have a unique optimal realization, whether such metrics 
can be recognized in polynomial time, and whether there exists a 
polynomial time algorithm for computing some optimal realization.  \\

\noindent{\bf Acknowledgment}: 
Hayamizu is supported by JST PRESTO Grant Number JPMJPR16EB. 
	% and is grateful to the organizers of Waiheke 2017.
	%, where she first introduced the concept of cactus metrics and $X$-cacti \mh{\cite{hay}}.} 
Huber, Moulton and Murakami thank the Netherlands Organization for 
Scientific Research (NWO), including Vidi grant 639.072.602. 
Huber and Moulton also thank the Research Institute for Mathematical Sciences, Kyoto University, 
	The Institute of Statistical Mathematics, Tokyo, and the London Mathematical Society for their support.

%\bigskip
%\noindent
%{\bf References}

%\bibliographystyle{elsarticle-num} 
%\bibliography{phylogeneticNetworks}

\begin{thebibliography}{10}
\expandafter\ifx\csname url\endcsname\relax
  \def\url#1{\texttt{#1}}\fi
\expandafter\ifx\csname urlprefix\endcsname\relax\def\urlprefix{URL }\fi
\expandafter\ifx\csname href\endcsname\relax
  \def\href#1#2{#2} \def\path#1{#1}\fi


\bibitem{alt}
I.Alth\"ofer, 
On optimal realizations of finite metric spaces by graphs, 
Discrete and Computational Geometry 3(1) (1988) 103-122.

\bibitem{rubei18}
A.Baldisserri, R.Elena,
Distance matrices of some positive-weighted graphs,
Australian Journal of Combinatorics  70(2) (2018) 185-201.

%\bibitem{DHK08b}
%	Dress AW, Huber KT, Koolen J, Moulton V. Cut points in metric spaces. 
% Applied Mathematics Letters. 2008 Jun 1;21(6):545-8.

\bibitem{band}
H.-J.Bandelt, Recognition of tree metrics, 
SIAM Journal on Discrete Mathematics 3(1) (1990) 1-6.

\bibitem{bat}
V.Batagelj, T.Pisanski, J.M.Sim\~{o}es-Pereira,
 An algorithm for tree-realizability of distance matrices, 
 International Journal of Computer Mathematics 34(3-4) (1990) 171-176.

\bibitem{culb}
J.C.Culberson, P.Rudnicki, A fast algorithm for constructing 
trees from distance matrices, Information Processing Letters 30(4) (1989) 215-220.

\bibitem{dress}
A.Dress,
Trees, tight extensions of metric spaces, and the cohomological dimension of certain groups: 
a note on combinatorial properties of metric spaces, 
Advances in Mathematics 53(3) (1984) 321-402.

\bibitem{dress10}
A.Dress, K.T.Huber, J.Koolen, V.Moulton, A.Spillner, 
An algorithm for computing cutpoints in finite metric spaces,
Journal of Classification 27(2) (2010) 158-172.

% \bibitem{gambette}
% Gambette P, Huber KT, Scholz GE. Uprooted Phylogenetic Networks. 
%	Bulletin of mathematical biology. 2017 Sep 1;79(9):2022-48

%\bibitem{gamb}
%P.Gambette, V.Berry, C.Paul, 
%Quartets and unrooted phylogenetic networks,
%Journal of Bioinformatics and Computational Biology   10(04) (2012) 1250004.

%\bibitem{geller}
%D.Geller, B.Manvel, 
%Reconstruction of cacti, 
%Canadian Journal of Mathematics  21 (1969) 1354-1360. 

\bibitem{hak}
S.L.Hakimi, S.S.Yau,
Distance matrix of a graph and its realizability, 
Quarterly of Applied Mathematics 22(4) (1965) 305-317.

\bibitem{hay}
M.Hayamizu, 
$X$-cactus trees and cactus tree metrics, The 21st New Zealand Phylogenomics Meeting (Waiheke 2017), 12-17 February 2017.
https://cdn.auckland.ac.nz/assets/compevol/events/\\documents/Waiheke2017Programme.pdf

%\bibitem{HVbridge}
%A.Hertz, S.Varone, 
%The metric bridge partition problem: partitioning of a metric space into two 
%subspaces linked by an edge in any optimal realization,
%Journal of Classification 24(2) (2007) 235-249.

\bibitem{HVcutpoint}
A.Hertz, S.Varone, 
The metric cutpoint partition problem,
Journal of Classification  25(2) (2008) 159-175.

%\bibitem{quarnet}
%K.T.Huber, V.Moulton, C.Semple, T.Wu, 
%Quarnet inference rules for level-1 networks, 
%Bulletin of Mathematical Biology 80(8) (2018) 2137-2153.

\bibitem{huson}
D.H.Huson, R.Rupp, C.Scornavacca, 
Phylogenetic networks: concepts, algorithms and applications, 
Cambridge University Press, 2010.

%	\bibitem{leo}
%	Van Iersel L, Moulton V. Leaf-reconstructibility of phylogenetic networks. 
%  SIAM Journal on Discrete Mathematics. 2018 Aug 7;32(3):2047-66

\bibitem{imrich84}
W.Imrich, J.M.Simoes-Pereira, C.M.Zamfirescu, 
On optimal embeddings of metrics in graphs,
Journal of Combinatorial Theory, Series B  36(1) (1984) 1-15.

%\bibitem{pbook}
%P.Lemey, M.Salemi, A.M.Vandamme, editors, 
%The phylogenetic handbook: a practical approach to phylogenetic analysis 
%and hypothesis testing, Cambridge University Press, 2009.

\bibitem{mat}
J.Matou{\v{s}}ek, 2.7 How large graph?, Open problems on embeddings of finite metric spaces
Workshop on discrete metric spaces and their algorithmic applications, 2002, available at http://kam.mff.cuni.cz/~matousek/metrop.ps. 

%\bibitem{nj}
%N. Saitou, M. Nei, 
%The neighbor-joining method a new method for reconstructing phylogenetic trees, 
%Molecular biology and evolution 4(4) (1987) 406-425. 

\bibitem{ss13}
C.Semple, M.Steel, Phylogenetics, Oxford University Press, 2003.

\bibitem{unicyclic}
J.M.S. Sim{\~o}es-Pereira, 
A note on distance matrices with unicyclic graph realizations,
Discrete Mathematics 65(3) (1987) 277-287.

\bibitem{compactification}
J.M.S. Sim{\~o}es-Pereira, C.M. Zamfirescu, 
Submatrices of non-tree-realizable distance matrices, 
Linear Algebra and its Applications 44 (1982) 1-17.

\bibitem{wink}
P.Winkler, 
The complexity of metric realization, 
SIAM Journal on Discrete Mathematics, 1(4) (1988) 552-559.

\end{thebibliography}

\end{document}